\documentclass[10pt]{amsart}

\title{On The Skorokhod Representation Theorem}
\author{Jean Cortissoz}

\newtheorem{theo}{Theorem}
\newtheorem{cor}{Corollary}
\newtheorem{theorem}{Theorem}[section]
\newtheorem{corollary}{Corollary}[section]
\newtheorem{lemma}{Lemma}[section]
\newtheorem{proposition}{Proposition}[section]
\newtheorem{definition}{Definition}[section]

\begin{document}
\begin{abstract}
In this paper we present a variant of the well known Skorokhod Representation Theorem. In our main result,
given $S$ a Polish Space,
to a given continuous path $\alpha$ in the space of probability measures on $S$, we associate a continuous path 
in the space of $S$-valued random variables on a nonatomic probability space (endowed with the topology of the convergence in probability).
We call this associated path a lifting of $\alpha$. 
An interesting feature of our result is that we can fix the endpoints of the lifting of $\alpha$, as long as their distribution
correpond to the respective endpoints of $\alpha$. We also discuss an $n$-dimensional generalization of this result.
\end{abstract}

\maketitle

\section{Introduction}
\label{introduction}

Let $\left(S,d\right)$ be a complete separable metric space and $\left(\Omega, \mathcal{F}, P\right)$ be a 
complete non atomic probability space (here $\mathcal{F}$ denotes the $\sigma$-algebra where $P$ is defined),
and $\mathcal{P}\left(S\right)$ the space of probability measures on $S$. The Skorokhod
Representation Theorem states the following
\begin{theo}
\label{Skorokhod}
Suppose $P_n$, $n=1, 2, \dots $ and $P$ are probability measures on $S$ (provided with
its Borel $\sigma$-algebra) such that
$P_n \Rightarrow P$ (see Section \ref{weakconvergence}, Definition \ref{wconvergence}). 
Then there is a probability space $\left(\Omega, \mathcal{F}, P\right)$ on which 
are defined $S$-valued random variables $X_n$, $n=1, 2, \dots$ and $X$ with distributions $P_n$ and $P$
respectively, such that $\lim_{n\rightarrow n}X_n =  X \quad a.s.$  
\end{theo}

A stronger result is presented in \cite{Blackwell}
The purpose of this paper is to prove a result in the same spirit of Theorem \ref{Skorokhod}. The main result we 
prove in this note is,

\begin{theo}
\label{maintheorem}
Let $\alpha:\,\left[0,1\right]\rightarrow \mathcal{P}\left(S\right)$ be a continuous function ($\mathcal{P}\left(S\right)$
endowed with the topology of the weak convergence -see Section \ref{weakconvergence}). Let $X_{\mu}$ and $X_{\nu}$ be random variables
such that $law\left(X_{\mu}\right)=\alpha\left(0\right)$ and $law\left(X_{\nu}\right)=\alpha\left(0\right)$.
There is $\hat{\alpha}:\,\left[0,1\right]\rightarrow L^0\left(\Omega,S\right)$ continuous ($L^0\left(\Omega,S\right)$ endowed with the
topology of the convergence in probability), such that $\hat{\alpha}\left(0\right)=X_{\mu}$ and $\hat{\alpha}\left(1\right)=X_{\nu}$
and $law\left(\hat{\alpha}\right)=\alpha$.
\end{theo}

A result with the same statement as our main theorem but with $\left[0,1\right]$ replaced by $\left[0,1\right]^n$ can be obtained
using techniques similar to the ones used in this paper. More exactly we have,
\begin{theo}
\label{maintheorem2}
Let $\alpha:\,\left[0,1\right]^n\rightarrow \mathcal{P}\left(S\right)$ be a continuous function ($\mathcal{P}\left(S\right)$
endowed with the topology of the weak convergence). Let $X:\,\partial \left[0,1\right]^n\rightarrow L^0\left(\Omega,S\right)$ 
($\partial$ = boundary)
be a continuous function
such that $law\left(X\right)=\alpha |_{\partial\left[0,1\right]^n}$. 
There is $\hat{\alpha}:\,\left[0,1\right]^n\rightarrow L^0\left(\Omega,S\right)$ continuous, such that $\hat{\alpha}|_{\partial \left[0,1\right]^n}=X$ 
and $law\left(\hat{\alpha}\right)=\alpha$.

\end{theo}

We must point out that the result in \cite{Blackwell} may seem stronger to the results described above. In certain sense this is
true, since given $\alpha:\,\left[0,1\right]\rightarrow \mathcal{P}\left(S\right)$ 
and a representation $\beta:\,\mathcal{P}\left(S\right)\times \Omega\longrightarrow S$ (using the definitions in \cite{Blackwell} -notice
that their $M\left(X\right)$ is our $\mathcal{P}\left(S\right)$)
by defining 
\[
\hat{\alpha}\left(t\right)=\beta\left(\alpha\left(t\right),\cdot\right)
\]

With this definition, we have that if $t_n\rightarrow t_0$ then $\hat{\alpha}\left(t_n\right)\rightarrow \hat{\alpha}\left(t_0\right)$ a.s.,
and a.s. convergence implies convergence in probability. However, in our results we are allowed to fix 
boundary values for the liftings (or representations). As a curious consequence of Theorem \ref{maintheorem2}  we obtain

\begin{cor}
All the homotopy groups of $L^0\left(\Omega,S\right)$ (endowed with the topology of convergence in probability) are trivial. 
\end{cor}

Even though we do not give a detailed argument to prove Theorem 
\ref{maintheorem2}, we present a rough sketch of a proof in the last sections
and details to be given in a further paper. This paper is arranged as follows: in Section \ref{BasicConcepts} we collect some 
basic theory as a quick reference for the convenience of the reader. In Section \ref{ImportantLemma} we prepare some lemmas we
use in the proof of Theorem \ref{maintheorem}, which
we finally prove in Section \ref{RepresentationTheorems}. The sketch of a proof of Theorem \ref{maintheorem2} at the end
of Section \ref{RepresentationTheorems}.

%%%%%%%%%%%%%%%%%%%%%%%%%%%%%%%%%%%%%%%%%%%%%%%%%%%%%%%%%%%%%%%%%%%%%%%%%%%%%%%%%%%%%%%%%%%%%%%%%%%%%%%%%%
\section{Some Basic Concepts}
\label{BasicConcepts}

In this section we collect some important and well known definitions and facts for the reader's convenience. We recommend Chapter 3 of
\cite{EthierKurtz} as a reference for this section.

\subsection{Convergence in Probability.}
\label{ProbabilityConvergence}
\begin{definition}

We say that a sequence $\left(X_n\right)$ converges in probability to $X$, and we denote it by
$X_n\stackrel{P}{\rightarrow} X$, if for every $\epsilon>0$ we have
\[
\lim_{n\rightarrow\infty} P\left\{\omega:\, d\left(X \left(\omega\right),Y\left(\omega\right)\right)\geq \epsilon\right\} = 0
\]

\end{definition}

\begin{definition}
Given $X$ and $Y$ random variables, define
\[
\rho\left(X,Y\right)=\inf\left\{\epsilon>0:\,
P\left\{\omega:\, d\left(X_n \left(\omega\right),Y\left(\omega\right)\right)\geq \epsilon\right\}\leq \epsilon\right\} 
\]
\end{definition}

\begin{theorem}
$\rho$ is a metric on $L^0\left(\Omega,S\right)$, and given a sequence $\left(X_n\right)$ 
of random variables and a random
variable $X$, then
\[
\lim_{n\rightarrow\infty}\rho\left(X_n,X\right)=0\quad\mbox{if and only if} \quad X_n \stackrel{P}{\rightarrow} X .
\]
\end{theorem}

\subsection{The Space $\mathcal{P}\left(S\right)$ and the function $law$. }
\label{weakconvergence}

\begin{definition}
\[
\mathcal{P}\left(S\right)=\left\{\mu:\, \mu \,\, \mbox{is a probability measure on}\,\, \mathcal{B}\left(S\right)\right\}.
\]
\end{definition}

\begin{definition}
Given $\mu \in\mathcal{P}\left(S\right)$ and $A\in\mathcal{B}\left(S\right)$ we say that $A$ is
a set of $\mu$-continuity if $\mu\left(\partial A\right)=0$ ($\partial A$ is the topological boundary of $A$).

A sequence $\left(\mu_n\right)_{n\in\mathbf{N}}$ in $\mathcal{P}\left(S\right)$ converges weakly to $\mu$,
and we denote it by $\mu_n\Rightarrow \mu$, if for every set of $\mu$-continuity $A$ 
$\mu_n\left(A\right)\rightarrow \mu\left(A\right)$ as $n\rightarrow \infty$.
\end{definition}

\begin{definition}
\label{wconvergence}

Given $\mu,\nu\in \mathcal{P}\left(S\right)$ define
\[
q\left(\mu,\nu\right)=\inf\left\{\epsilon>0:\, \mu\left(A\right)\leq \nu\left(A^{\epsilon}\right)+\epsilon
\quad \mbox{for all} \quad A\subset S \,\, \mbox{closed}\right\}
\]
\end{definition}

\begin{theorem}
$q$ defines a metric on $\mathcal{P}\left(S\right)$, and given a sequence $\left(\mu_n\right)$ and
a probability measure $\mu$, then
\[
\lim_{n\rightarrow \infty}q\left(\mu_n,\mu\right)=0 \quad \mbox{if and only if} \quad
\mu_n\Rightarrow \mu.
\]
\end{theorem}

\begin{definition}
Given $X\in L^0\left(\Omega,S\right)$, the probability measure $\mu$ defined on $\mathcal{B}\left(S\right)$ by
\[
\mu\left(A\right)=P\left(X^{-1}\left(A\right)\right)
\]
is called the {\bf distribution} or {\bf law} of $X$ .

\end{definition}

%%%%%%%%%%%%%%%%%%%%%%%%%%%%%%%%%%%%%%%%%%%%%%%%%%%%%%%%%%%%%%%%%%%%%%%%%%%%%%%%%%%%%%%%%%%%%%%%%%%%%%%%%%%%
\section{An important lemma}
\label{ImportantLemma}

\begin{theorem}
\label{strassen1}
Let $\left(S,d\right)$ be separable, and let $P,Q \in \mathcal{P}\left(M\right)$.
Define $\mathcal{M}\left(P,Q\right)$ be the set of all $\mathcal{P}\left(S\times S\right)$
with marginals $P$ and $Q$. Then
\[
q\left(P,Q\right)=\inf_{\mu\in \mathcal{M}\left(P,Q\right)}
\inf\left\{\epsilon>0:\, \mu\left\{\left(x,y\right):\,d\left(x,y\right)\geq \epsilon\right\}\leq \epsilon\right\}
\]
\end{theorem}

\bfseries\textit{Proof. }\normalfont  See \cite{EthierKurtz} (Chapter 3, Theorem 1.2) or \cite{Strassen} (Corollary to Theorem 10).

\hfill $\Box$

As a corollary we get

\begin{corollary}
\label{strassen2}
\[
q\left(P,Q\right)=\inf\left\{\rho\left(X,Y\right):\, law\left(X\right)=P \quad \mbox{and} \quad law\left(Y\right)=Q \right\} .
\]
\end{corollary}

\bfseries\textit{Proof. }\normalfont It follows from the fact that given $\mu\in \mathcal{M}\left(P,Q\right)$ there is a 
random variable $W\in L^0\left(\Omega,S\times S\right)$ such that $law\left(W\right)=\mu$. But by its very definition
$W=\left(X,Y\right)$ where $X,Y\in L^0\left(\Omega,S\right)$ and $law\left(X\right)=P$ and $law\left(Y\right)=Q$.
Then all we have to notice is that
\[
\rho\left(X,Y\right)
=\inf \left\{\epsilon>0:\, P\left\{\omega:\,d\left(X\left(\omega\right),Y\left(\omega\right)\right)\geq\epsilon \right\}\leq \epsilon \right\} .
\]

\hfill $\Box$

We will need the following ``working'' definitions, 

\begin{definition}
$\delta_a$ is the measure defined by 
\[
\delta_a \left(V\right)=\left\{
\begin{array}{l}
1 \quad \mbox{if}\quad a\in V\\
0 \quad \mbox{otherwise}
\end{array}
\right.
\]

We say that $\mu$ is finitely supported if it can be written as
\[
\mu = \sum_{i=1}^n c_i \delta_{a_i} \quad c_i\geq 0
\]
\end{definition}

\begin{definition}
Let $\left(A_k\right)_{k=1,\dots,n}$ be a partition of  $\Omega$.
We define the simple $S$-valued random variable $X=\sum_{i=1}^n \chi_{A_i}^{a_i}$ as
\[
X\left(\omega\right)=a_i \quad \mbox{iff} \quad \omega\in A_i
\]

\end{definition}

Corollary \ref{strassen2} is used to prove

\begin{lemma}
\label{strassen3}
Let $\mu$ and $\nu$ be finitely supported measures, and let $\epsilon>0$ be such that
\[
q\left(\mu,\nu\right)<\epsilon .
\]

Then, given X a random variable such that $law\left(X\right)=\mu$ there exists $Y$ such that
$law\left(Y\right)=\nu$ and $\rho\left(X,Y\right)<\epsilon$
\end{lemma}

\bfseries\textit{Proof. }\normalfont 
Since $\mu$ and $\nu$ have finite support, by Theorem \ref{strassen2} we can find 
\[
X'=\sum_{j=1}^m \chi_{A_j'}^{a_j} \quad \mbox{and} \quad
Y'=\sum_{j=1}^m \chi_{B_j'}^{a_j} .
\]
so that $law\left(X'\right)=\mu$ and $law\left(Y'\right)=\nu$
and $\rho\left(X',Y'\right)<\epsilon$. 

Write $X=\sum_{j=1}^m \chi_{A_j}^{a_j}$. Since $\Omega$ is nonatomic we can find measurable sets
$B_1,B_2,\dots, B_m$ such that
\[
P\left(A_i\cap B_j\right)=P\left(A_i'\cap B_j'\right) \quad \mbox{for all} \quad i,j=1, 2,\dots, m .
\]

It is clear then that $Y=\sum_{j=1}^m \chi_{B_j}^{a_j}$ satisfies $\rho\left(X,Y\right)<\epsilon$.

\hfill $\Box$

Finally we have the following fundamental lemma,
 
\begin{lemma}
Let $\epsilon>0$ be given and assume $q\left(law\left(X\right),law\left(Y\right)\right)<\epsilon$. Then there
is $Y'$ such that $law\left(Y'\right)=law\left(Y\right)$ and $\rho\left(X,Z\right)<\epsilon$.
\end{lemma}

\bfseries\textit{Proof. }\normalfont
We proceed by induction on the complexity of random variables. The case when $X$ and
$Y$ are simple is contained in corollary \ref{strassen2}.

Now assume that $X$ is arbitrary and $Y$ simple. Let $\delta>0$ be such that $\rho\left(X,Y\right)<\epsilon-\delta$. 
Since the set of simple random variables is dense in $L^0\left(\Omega,S\right)$,
we can choose $X'$ a simple random variable such that $\rho\left(X,X'\right)<\delta$. By the induction hypothesis,
we can find $Y'$ such that $law\left(Y'\right)=law\left(Y\right)$ and $\rho\left(X',Y'\right)<\epsilon-\delta$.
Hence we have,

\[
\rho\left(X,Y'\right)\leq \rho\left(X,X'\right)+\rho\left(X',Y'\right)<\delta+\epsilon-\delta=\epsilon
\] 

Finally, let $X$ and $Y$ be arbitrary random variables and let $\delta>0$ be such that

\[
q\left(law\left(X\right),law\left(Y\right)\right)+1000\delta<\epsilon .
\]
 
Choose a sequence $\left(Y_n\right)$ of simple random
variables converging to $Y$ and such that 
\[
\rho\left(Y_n,Y_{n+1}\right)<\frac{1}{2^{n+1}}\quad \mbox{and}\quad 
q\left(law\left(Y_n\right),law\left(Y\right)\right)<\delta ,
\]
and let $N$ be such that $\frac{1}{2^N}<\delta$ and also $\rho\left(Y_N,Y\right)<\delta$.
We construct a new sequence $\left(Y_j'\right)_{j=N,N+1,N+2,\dots}$ as follows:

If $j=N$, choose $Y_N'$ be such that $law\left(Y_N'\right)=law\left(Y_N\right)$ and
\[
\rho\left(Y_N',X\right)<q\left(law\left(X\right),law\left(Y\right)\right)+\delta
\]

This can be done because
\[
\begin{array}{rcl}
q\left(law\left(Y_N\right),law\left(X\right)\right)&\leq& q\left(law\left(Y_N\right),law\left(Y\right)\right)
+q\left(law\left(Y\right),law\left(X\right)\right)\\
&\leq& q\left(law\left(Y\right),law\left(X\right)\right)+\delta .
\end{array}
\]

Once we have chosen $Y_j'$ for $j=N+1,\dots,N+M$, we pick $Y_{N+M+1}'$ such that 
\[
law\left(Y_{N+M+1}\right)=law\left(Y_{N+M+1}'\right) \quad
\mbox{and} \quad \rho\left(Y_{N+M+1}',Y_{N+M}'\right)<\frac{1}{2^{N+M+1}}
\]

This can be done because
\[
q\left(law\left(Y_{N+M+1}\right),law\left(Y_{N+M}'\right)\right)=
q\left(law\left(Y_{N+M+1}\right),law\left(Y_{N+M}\right)\right)\leq \frac{1}{2^{N+M+2}} .
\]

By construction, the sequence $\left(Y_j'\right)_j$ is convergent, and for its limit $Y'$
it holds that $law\left(Y'\right)=law\left(Y\right)$, and also
\[
\begin{array}{rcl}
\rho\left(X,Y'\right)&\leq& \rho\left(X,Y_N'\right)+\sum_{j=1}\frac{\delta}{2^j}\\
&<& q\left(law\left(X\right),law\left(Y\right)\right)+\delta+\delta<\epsilon .
\end{array}
\]

\hfill $\Box$

\section{Representation Theorems}
\label{RepresentationTheorems}

\subsection{Liftings. } 
\label{Liftings}

First we give a definition we learnt from Ramiro de la Vega.
\begin{definition}
A family of measurable sets $\left(A_t\right)_{t\in\left[0,\delta\right]}$ is a $\left[0,\delta\right]$-family
if it satisfies:

(i) $A_s\subset A_t$ whenever $s\leq t$, 

(ii) $P\left(A_t\right)=t$ .
\end{definition}

The following lemma we also learnt from de La Vega, is what makes $\left[0,\delta\right]$-families a useful tool.

\begin{lemma}
Let $\left(\Omega,\mathcal{F}, P\right)$ be a complete nonatomic probability space, and let $A\in \mathcal{F}$ and
let $\delta=P\left(A\right)$. Then there is a $\left[0,\delta\right]$-family $\left(A_t\right)_{t\in\left[0,\delta\right]}$
such that $A_t\subset A$ for all $t$.
\end{lemma}

\bfseries\textit{Proof. }\normalfont Pick an ordering of $\mathbf{Q}\cap \left[0,1\right]$ ($\mathbf{Q}$: the rational numbers) say 
$q_1, q_2, q_3, \dots$ . Since $\Omega$ is nonatomic we can inductively construct $E_{q_1}, E_{q_2}, E_{q_3}, \dots$
such that 
\[
P\left(E_{q}\right)=q, \quad \mbox{and}\quad E_{q}\subset E_{r}\quad \mbox{if}\quad q\leq r . 
\]

Finally for $x\in \left[0,1\right]\setminus \mathbf{Q}$, define
\[
E_x = \bigcup_{q<x}E_q
\]
\hfill $\Box$

Using $\left[0,\delta\right]$ families we can introduce the notion of a segment joining two simple random variables.

\begin{definition}
\label{segments1}
Let
\[
X=\sum_{j=1}^m \chi_{A_j}^{a_j} \quad \mbox{and} \quad
Y=\sum_{j=1}^m \chi_{B_j}^{a_j}
\]
be two simple random variables. Define
\[
E_{ij}=A_i\cap B_j \quad\mbox{and} \quad e_{ij}=P\left(E_{ij}\right) .
\]
and let $\left(\left[E_{ij}\right]_t\right)$ be a $\left[0,e_{ij}\right]$ family of $E_{ij}$. A segment 
$\alpha_{X,Y}:\,\left[a,b\right]\rightarrow L^0\left(\Omega,S\right)$
joining $X$ and $Y$ is defined as

\[
\alpha_{X,Y}\left(t\right)= \sum_{i=1}^m \chi^{a_i}_{E_{ii}\cup\left(\bigcup_{k=1,k\neq i} 
\left[E_{ki}\right]_{\left(\frac{t-a}{b-a}\right)e_{ki}}\right)}
+\sum_{i=1}^m\sum_{j=1, j\neq i}^m \chi_{E_{ij}\setminus \left[E_{ij}\right]_{\left(\frac{t-a}{b-a}\right)e_{ij}}}^{a_i} 
\]

\end{definition} 

We describe some important properties of these segments.

\begin{proposition}
\label{liftingsegments}
$\hat{\alpha}:=\alpha_{X,Y}$ thus defined is a continuous function with $\alpha_{X,Y}\left(a\right)=X$ and 
$\alpha_{X,Y}\left(b\right)=Y$. Moreover,
\[
\alpha\left(t\right):=law\left(\hat{\alpha}\left(t\right)\right)=\left(\frac{b-t}{b-a}\right)law\left(X\right)
+\left(\frac{t-a}{b-a}\right)law\left(Y\right) .
\] 
\end{proposition}

\bfseries\textit{Proof. }\normalfont First we show that $\alpha_{X,Y}$ is continuous. It is an
immediate consequence of the following inequality. Let $\epsilon>0$ be given and $s\leq t$, then we have,
\[
\begin{array}{rcl}
P\left\{\omega:\, d\left(\hat{\alpha}\left(t\right),\hat{\alpha}\left(s\right)\right)\geq \epsilon\right\}&\leq&
\sum_{\left\{\left(i,k\right):\, d\left(a_i,a_k\right)\geq \epsilon\right\}} P\left(\left[E_{ik}\right]_{\left(t-s\right)e_{ik}}\right)\\
&\leq& \left(t-s\right)\sum_{\dots} e_{ik}\leq t-s .
\end{array}
\]

Let us show that $law\left(\hat{\alpha}\right)=\alpha$. To make things easier, we will assume $a=0$ and $b=1$.
Then all we must show is that the coefficient of $\delta_{a_i}$ in $law\left(\hat{\alpha}\right)$ is
$\left(1-t\right)P\left(A_i\right)+tP\left(B_i\right)$. Let's fix $i=1$. Then the sought coefficient is given by
\[
\begin{array}{c}
P\left(E_{11}\right)+\sum_{j=2}^m P\left(\left[E_{j1}\right]_{te_{j1}}\right)+
\sum_{j=2}^m P\left(E_{1j}\setminus\left[E_{1j}\right]_{tm_{1j}}\right)\\
=e_{11}+\sum_{j=2}^m tm_{j1}+\sum_{j=2}^m \left(m_{1j}-tm_{1j}\right)\\
=P\left(A_1\right)+t\left(P\left(B_1\right)-e_{11}\right)-t\sum_{j=2}^m e_{1j}\\
=P\left(A_1\right)+t\left(P\left(B_1\right)-e_{11}\right)-t\sum_{j=2}^m e_{1j}\\
=P\left(A_1\right)+tP\left(B_1\right)-t\sum_{j=1}^m e_{1j}\\
=P\left(A_1\right)+tP\left(B_1\right)-tP\left(A_1\right) .
\end{array}
\]

\hfill $\Box$

The notion of segments can be generalized to the concept of a ``poligonal''.
\begin{definition}
We call $\beta:\, \left[0,1\right]\longrightarrow \mathcal{P}\left(S\right)$ a poligonal with
vertices at $\mu_0, \mu_1,\dots,\mu_n,\mu_{n+1}$ if there is a partition $0=t_0,t_1,\dots, t_n,t_{n+1}=1$
of $\left[0,1\right]$ such that $\beta$ restricted to $\left[t_i,t_{i+1}\right]$ is given by
\[
\beta\left(t\right)=\left(\frac{t_{i+1}-t}{t_{i+1}-t_i}\right)\mu_i+\left(\frac{t-t_i}{t_{i+1}-t_i}\right)\mu_{i+1} .
\]
\end{definition}

An easy consequence of proposition \ref{liftingsegments} is the following fact about poligonals,

\begin{proposition}
\label{liftingpolygonals0}
Let $\alpha:\,\left[0,1\right]\longrightarrow \mathcal{P}\left(S\right)$ be a poligonal with 
vertices at measures of finite support, and let $\alpha\left(0\right)=\mu$ and $\alpha\left(1\right)=\nu$.
Given $X_{\mu}$ and $X_{\nu}$ such that $law\left(X_{\mu}\right)=\mu$ and $law\left(X_{\nu}\right)=\nu$,
there is a lifting $\hat{\alpha}:\,\left[0,1\right]\longrightarrow L^0\left(\Omega,S\right)$ (i.e., $law\left(\hat{\alpha}\right)=\alpha$),
such that $\hat{\alpha}\left(0\right)=X_{\mu}$ and $\hat{\alpha}\left(1\right)=X_{\nu}$.

\end{proposition}

Also poligonals are dense in the space of continuous maps from the unit interval to $\mathcal{P}\left(S\right)$.
Before we write and prove the exact statement of this fact we will need the following observation.

\begin{lemma}
\label{q-distance}
Let $\mu,\nu \in \mathcal{P}\left(S\right)$. For $t\in\left[0,1\right]$ we have
\[
q\left(\nu, t\mu+\left(1-t\right)\nu\right)\leq q\left(\nu,\mu\right).
\] 
\end{lemma}

\bfseries\textit{Proof. }\normalfont Let $\epsilon>0$ be such that $\mu\left(A\right)\leq \nu\left(A^{\epsilon}\right)+\epsilon$
for all $A\subset S$ closed. Then we have,
\[
\begin{array}{rcl}
t\mu\left(A\right)+\left(1-t\right)\nu\left(A\right)&\leq& t\nu\left(A^{\epsilon}\right)+t\epsilon+
\left(1-t\right)\nu\left(A^{\epsilon}\right)+\left(1-\epsilon\right)\epsilon\\
&=& \nu\left(A^{\epsilon}\right)+\epsilon
\end{array}
\]
and from this the statement of the lemma follows.

\hfill $\Box$

Now we are ready to state and prove the following density property of poligonals

\begin{lemma}
Given $\alpha:\,\left[0,1\right]\longrightarrow \mathcal{P}\left(S\right)$ and $\epsilon>0$ there is a poligonal
$\beta$ with vertices at measures of finite support such that
\[
\sup_{t\in\left[0,1\right]}q\left(\alpha\left(t\right),\beta\left(t\right)\right)<\epsilon
\]
\end{lemma}

\bfseries\textit{Proof. }\normalfont Let $\epsilon>0$ be given. By the uniform continuity of $\alpha$, we can find $\delta>0$
such that whenever $\left|s-t\right|<\delta$ we have $q\left(\alpha\left(t\right),\alpha\left(s\right)\right)<\frac{\epsilon}{5}$.
Let $N>0$ be big enough so that $\frac{1}{N}<\delta$, and define a partition of the interval $\left[0,1\right]$ by
$t_i=\frac{i}{N}$ $i=0,1,\dots,N$. For each $i$ pick a finitely supported measure $\mu_i$ such that 
$q\left(\mu_i,\alpha\left(t_i\right)\right)\leq \frac{\epsilon}{5}$. Let $\beta$ be the poligonal defined by the segments
$\beta:\,\left[t_i,t_{i+1}\right]\longrightarrow \mathcal{P}\left(S\right)$ with endpoints $\mu_i$ and $\mu_{i+1}$.
For $t\in \left[t_i,t_{i+1}\right]$ we have,
\[
\begin{array}{rcl}
q\left(\alpha\left(t\right),\beta\left(t\right)\right)&\leq& q\left(\alpha\left(t\right),\alpha\left(t_i\right)\right)
+q\left(\alpha\left(t_i\right),\mu_i\right)+q\left(\mu_i,\beta\left(t\right)\right)\\
&\mbox{by Lemma \ref{q-distance}}& \\
&\leq& q\left(\alpha\left(t\right),\alpha\left(t_i\right)\right)
+q\left(\alpha\left(t_i\right),\mu_i\right)+q\left(\mu_i,\mu_{i+1}\right)\\
&\leq& q\left(\alpha\left(t\right),\alpha\left(t_i\right)\right)
+q\left(\alpha\left(t_i\right),\mu_i\right)+q\left(\mu_i,\alpha\left(t_i\right)\right)\\
&& + q\left(\alpha\left(t_i\right),\alpha\left(t_{i+1}\right)\right)
+q\left(\alpha\left(t_{i+1}\right),\mu_{i+1}\right)\\
&\leq& \frac{\epsilon}{5}+\frac{\epsilon}{5}+\frac{\epsilon}{5}+\frac{\epsilon}{5}+\frac{\epsilon}{5}+\frac{\epsilon}{5}=\epsilon .
\end{array}
\]

\hfill $\Box$

\subsection{Proof of the Main Theorem.} We are almost ready to prove the Main Theorem of this paper. The following
fact will be used in its proof.

\begin{lemma}
\label{liftingpoligonals}
Let $\alpha:\, \left[0,1\right]\rightarrow \mathcal{P}\left(S\right)$ and  let $\epsilon>0$ be given. Let $\beta$ be an arbitrary poligonal
with vertices at measures of finite support and such that
\[
\sup_{t\in\left[0,1\right]}q\left(\alpha\left(t\right),\beta\left(t\right)\right)<\epsilon
\]

Then, given any continuous lifting $\hat{\alpha}$ of $\alpha$, there is a lifting $\hat{\beta}$ of $\beta$ such that
\[
\sup_{t\in\left[0,1\right]}\rho\left(\hat{\alpha}\left(t\right),\hat{\beta}\left(t\right)\right)< 5\epsilon .
\]
\end{lemma}

For the proof of this lemma we need the following observation,
\begin{lemma} 
\label{d-distance}
Let 
\[
X_{\mu}=\sum_{j=1}^m \chi_{A_j}^{a_j} \quad \mbox{and} \quad X_{\nu}=\sum_{j=1}^m \chi_{B_j}^{a_j}
\]
be such that $law\left(X_{\mu}\right)=\mu$ and $law\left(X_{\nu}\right)=\nu$
are finitely supported measures. Let 
$\hat{\alpha}:=\alpha_{X_{\mu},X_{\nu}}$ be as in Definition \ref{segments1}. Then we have
\[
\rho\left(X_{\mu},\hat{\alpha}\left(t\right)\right)\leq \rho\left(X_{\mu},X_{\nu}\right) .
\]
\end{lemma}

\bfseries\textit{Proof. }\normalfont (We use the notation of Definition \ref{segments1}) Given $\epsilon>0$ we have
\[
\begin{array}{rcl}
P\left\{\omega:\, d\left(X_{\mu}\left(\omega\right),X_{\nu}\left(\omega\right)\right)\geq \epsilon\right\}&=& 
\sum_{\left\{\left(i,j\right):\, d\left(a_i,a_j\right)\geq \epsilon \right\}}P\left(E_{ij}\right)\\
&\geq& 
\sum_{\left\{\left(i,j\right):\, d\left(a_i,a_j\right)\geq \epsilon \right\}} P\left(\left[E_{ij}\right]_{te_{ij}}\right)\\
&=& P\left\{\omega:\, d\left(X_{\mu}\left(\omega,\hat{\alpha}\left(t\right)\right)\right)\geq \epsilon\right\} .
\end{array}
\]

The conclusion of the lemma follows.

\hfill $\Box$

\bfseries \textit{Proof of Lemma \ref{liftingpoligonals}. }\normalfont
Let $\hat{\alpha}$ be a continuous lifting of $\alpha$.
Take a partition $0=t_0<t_1<\dots<t_{n+1}=1$ of the unit interval, in such a way that
\[
\rho\left(X_i,X_{i+1}\right)<\epsilon \quad \mbox{where}\quad X_i=\hat{\alpha}\left(t_i\right)
\]

Choose $Y_i$ for $i=0,1,\dots,n+1$ so that $law\left(Y_i\right)=\beta\left(t_i\right)$ and
$\rho\left(X_i,Y_i\right)<\epsilon$. Then we have
\[
\rho\left(Y_i,Y_{i+1}\right)\leq \rho\left(Y_i,X_i\right)+\rho\left(X_i,X_{i+1}\right)
+\rho\left(X_{i+1},Y_{i+1}\right)< 3\epsilon .
\]

Construct a lifting $\hat{\beta}$ of $\beta$, such that 
$\hat{\beta}$ restricted to the segment $\left[t_i,t_{i+1}\right]$ is a lifting of
$\beta:\,\left[t_i,t_{i+1}\right]\longrightarrow \mathcal{P}\left(S\right)$ with $\hat{\beta}\left(t_i\right)=Y_i$
as given by Definition \ref{segments1}. Then $\hat{\beta}$ is continuous and for $t_i\leq t<t_{i+1}$ we have
\[
\begin{array}{rcl}
\rho\left(\hat{\alpha}\left(t\right),\hat{\beta}\left(t\right)\right)&\leq&
\rho\left(\hat{\alpha}\left(t\right),X_i\right)+\rho\left(X_i,Y_i\right)+ \rho\left(Y_i,\hat{\beta}\left(t\right)\right)\\
&\mbox{by Lemma \ref{d-distance}} & \\
&\leq& \rho\left(\hat{\alpha}\left(t\right),X_i\right)+\rho\left(X_i,Y_i\right)+ \rho\left(Y_i,Y_{i+1}\right)\\
&\leq& \epsilon+\epsilon+3\epsilon = 5\epsilon .
\end{array}
\]

\hfill $\Box$

\bfseries\textit{Proof of Theorem \ref{maintheorem}. }\normalfont 
Take a sequence of poligonals $\left(\alpha_n\right)_{n\in \mathbf{N}}$ with vertices at measures of finite support, and such that
\[
\alpha_n\rightarrow \alpha \quad \mbox{and} \quad 
\sup_{t\in\left[0,1\right]}\left(q\left(\alpha_n\left(t\right),\alpha_{n+1}\left(t\right)\right)\right)<\frac{1}{5^{n+1}} .
\] 

By proposition
\ref{liftingpolygonals0} and lemma \ref{liftingpoligonals},
we can lift this sequence to a sequence $\left(\hat{\alpha}_n\right)$ so that $\hat{\alpha}_n\left(0\right)\rightarrow X_{\mu}$,
 $\hat{\alpha}_n\left(1\right)\rightarrow X_{\nu}$ and
\[
\sup\left(\rho\left(\alpha_n\left(t\right),\alpha_{n+1}\left(t\right)\right)\right)<\frac{1}{5^n}  .
\]

It is clear by construction that $\left(\hat{\alpha}_n\right)$ is a convergent sequence. Let $\hat{\alpha}$ be its limit. Then,
since the convergence is uniform, $\hat{\alpha}$ is continuous, and because $law$ is a continuous function, 
$law\left(\hat{\alpha}\right)=\alpha$. This finishes the proof.
 
\hfill $\Box$

\subsection{On Theorem \ref{maintheorem2}. } Here we say a couple of words on how to approach a proof for Theorem \ref{maintheorem2}, with
further details to be given in a further paper. First, we use a special family of functions to approximate continous maps.
Let
\[
g_{\left(\mu_1,\mu_2,\dots,\mu_n\right)}\,:\, \left[0,1\right]^n \longrightarrow \mathcal{P}\left(S\right)
\]
as follows. First we define,
\[
g_{\left(\mu_1,\mu_2\right)}=\left(1-t_1\right)\mu_1+t_1\mu_2
\]
and then inductively
\[
g_{\left(\mu_1,\mu_2,\dots,\mu_n,\mu_{n+1}\right)}=(1-t_{n+1})g_{\left(\mu_1,\mu_2,\dots,\mu_n\right)}+
t_{n+1}\mu_{n+1}
\]

We use the functions $g_{\left(\mu_1,\mu_2,\dots,\mu_n\right)}$ 
(or natural variations of them) to approximate continuous functions from $\left[0,1\right]^n$
to $\mathcal{P}\left(S\right)$ in the same way we use polygonals to approximate continuous functions from $\left[0,1\right]$
to $\mathcal{P}\left(S\right)$.
Therefore we must learn how to lift these functions. We use the following procedure inductive procedure. Let
\[
\hat{g}_{\left(\mu_1,\dots,\mu_n\right)}\left(t_1,\dots,t_n\right)=\sum \chi^{a_i}_{A_i\left(t_1,
\dots,t_n\right)}
\]
be a lift of $g_{\left(\mu_1,\dots,\mu_n\right)}$. Given $X_{n+1}=\sum \chi_{B_i}^{a_i}$ a lift
of $\mu_{n+1}$, we construct a lifting of $g_{\left(\mu_1,\dots,\mu_n,\mu_{n+1}\right)}$
as follows (we adopt the notation $\stackrel{\rightarrow}{t}:=\left(t_1,\dots,t_n\right)$, and
$\stackrel{\rightarrow}{\mu}:=\left(\mu_1,\dots,\mu_n\right)$). First for each $j$ choose a $\left[0,m_j\right]$-family
of $B_j$ (here $m_j=P\left(B_j\right)$).
For fixed $\stackrel{\rightarrow}{t}$ construct a $[0,e_{ij}\left(\stackrel{\rightarrow}{t}\right)]$-family of $E_{ij}\left(\stackrel{\rightarrow}{t}\right)
=A_i\left(\stackrel{\rightarrow}{t}\right)\cap B_j$, by taking
\[
\left[E_{ij}\left(\stackrel{\rightarrow}{t}\right)\right]_{\gamma}
= E_{ij}\left(\stackrel{\rightarrow}{t}\right)\cap 
\left[B_j\right]_{\sup\left\{s:\,P\left(E_{ij}\left(\stackrel{\rightarrow}{t}\right)\cap\left[B_j\right]_s\right)=\gamma\right\}} .
\] 

Then define,
\[
\begin{array}{rcl}
\hat{g}_{\left(\stackrel{\rightarrow}\mu,\mu_{n+1}\right)}\left(\stackrel{\rightarrow}{t},t_{n+1}\right)&=&\sum_i
\chi_{E_{ii}\left(\stackrel{\rightarrow}{t}\right)\cup\left(\bigcup_{k=1,k\neq i}
\left[E_{ki}\right]_{t_{n+1}\cdot e_{ki}\left(\stackrel{\rightarrow}{t}\right)}\right)}^{a_i}
\\
&&+\sum_{i=1}^m\sum_{j=1, j\neq i}^m 
\chi_{E_{ij}\left(\stackrel{\rightarrow}{t}\right)\setminus 
\left[E_{ij}\right]_{t_{n+1}\cdot e_{ij}\left(\stackrel{\rightarrow}{t}\right)}}^{a_i} .
\end{array}
\]

Of course we must show that the lifting thus defined is continuous. Once we have proved that this lifting is continuous, we approximate the given
function with prescribed boundary values using functions in this family of liftings by choosing the vertices of the
approximation wisely, and then taking a uniform limit of the approximations, as it was done for the case of the unit interval.

\end{document}